\documentstyle[12pt]{article}
\title{The First Zagreb Index Conditions for Some Hamiltonian Properties of Graphs}
\author {Rao Li\\
              Dept. of Computer Science, Engineering and Mathematics\\
              University of South Carolina Aiken\\
              Aiken, SC 29801 \\ 
              USA \\}
\date{Aug. 31, 2024}

\begin{document}
\maketitle
\begin{abstract}
Let $G = (V, E)$ be a graph. The first Zagreb index of a graph $G$ is defined as
$\sum_{u \in V} d^2(u)$, where $d(u)$ is  the degree of vertex $u$ in $G$.  Using the P\'{o}lya-Szeg\H{o} inequality, 
we in this paper present the first Zagreb index conditions for some Hamiltonian properties of a graph and an upper bound for the first Zagreb index of a graph.
\end{abstract}  
$$Mathematics \,\, Subject \,\, Classification: 05C45, \,\, 05C09$$
$$Keywords:   The \,\, first \,\, Zagreb \,\, index, \,\,  Hamiltonian \,\, graph, $$
$$traceable \,\, graph $$

 \noindent {\bf 1 \, Introduction } \\
 
\noindent We consider only finite undirected graphs without loops or multiple edges.
Notation and terminology not defined here follow those in \cite{Bondy}.
Let $G = (V(G), E(G))$ be a graph with $n$ vertices and $e$ edges, the degree of a vertex $v$ is denoted by $d_G(v)$. We use $\delta$ and $\Delta$ to denote the minimum degree and maximum degree of $G$, respectively.  
A set of vertices in a graph $G$ is independent if the vertices in the set are pairwise nonadjacent. A maximum independent set in a graph $G$ is an independent set of largest possible size. 
 The independence number, denoted $\beta(G)$, of a graph $G$ is the cardinality of a maximum independent set in $G$. 
 For disjoint vertex subsets $X$ and $Y$ of $V(G)$, we use 
 $E(X, Y)$ to denote the set of all the edges in $E(G)$ such that one end vertex of each edge is in $X$ and another end vertex of the edge is in $Y$. Namely, 
$E(X, Y) := \{\, e : e = xy \in E, x \in X, y \in Y \,\}$. 
A cycle $C$  in a graph $G$ is called a Hamiltonian cycle of $G$ if $C$ contains all the vertices of $G$. 
A graph $G$ is called Hamiltonian if $G$ has a Hamiltonian cycle.
A path $P$  in a graph $G$ is called a Hamiltonian path of $G$ if $P$ contains all the vertices of $G$. 
A graph $G$ is called traceable if $G$ has a Hamiltonian path.\\

The first Zagreb index of a graph was introduced by Gutman and Trinajsti\'{c} in \cite{GT}. For a graph $G$, its first Zagreb index,
denoted $M_1(G)$,
 is defined as $\sum_{u \in V(G)} d_G^{2}(u)$. 
  Using the P\'{o}lya-Szeg\H{o} inequality, we in this paper to present the first Zagreb index conditions  
for the Hamiltonian and traceable graphs and an upper bound for the first Zagreb index of a graph. The main results are as follows. \\

\noindent {\bf Theorem $1$.} Let $G$ be a $k$-connected ($k \geq 2$) graph with $n \geq 3$ vertices and $e$ edges.  
 If 
 $$M_1(G) \geq  (n - k - 1)\Delta^2 + \frac{(e (\delta + n - k - 1))^2}{4 \delta (n - k - 1)(k + 1)},$$
 then $G$ is Hamiltonian or $G$ is $K_{k, \, k + 1}$. \\ 
 
 \noindent {\bf Theorem $2$.} Let $G$ be a $k$-connected ($k \geq 1$) with $n \geq 9$ vertices and $e$ edges. If
  $$M_1(G) \geq  (n - k - 2)\Delta^2 + \frac{(e(\delta + n - k - 2))^2}{4 \delta (n - k - 2)(k + 2)},$$ 
 then $G$ is traceable or $G$ is $K_{k, \, k + 2}$. \\
 
 \noindent {\bf Theorem $3$.} Let $G$ be a graph with $n$ vertices, $e$ edges, and $\delta \geq 1$. Then \\
 $$M_1(G) \leq  (n - \beta)\Delta^2 + \frac{(e (\delta + n - \beta))^2}{4 \delta (n - \beta) \beta}$$
 with equality if and only if $G$ is $K_{\beta, \, n - \beta}$ or $G$ is a bipartite graph with partition sets of $I$ and $V - I$ such that $|I| = \beta$, $\delta <  n - \beta$, $d(v) = \Delta$ for each vertex $v$ in $v - I$, and 
$I = P \cup Q$, where $P = \{\, x : x \in I, d(x) = n - \beta \,\}$,  $Q = \{\, y : y \in I, d(y) = \delta \,\}$, 
$|P| = \frac{\delta \beta}{\delta + n - \beta}$ which is an integer, and $|Q| = \frac{(n - \beta)\beta}{\delta + n - \delta}$ which is is an integer. \\

\noindent {\bf 2 \, Lemmas } \\

We will use the following results as our lemmas. The first two are from \cite{CE}. \\

\noindent {\bf Lemma $1$} \cite{CE}. Let $G$ be a $k$-connected graph of order $n \geq 3$. If $\beta \leq k$, then $G$ is Hamiltonian. \\

\noindent {\bf Lemma $2$} \cite{CE}. Let $G$ be a $k$-connected graph of order n. If $\beta \leq k + 1$, then $G$ is traceable. \\

Lemma $3$ is the P\'{o}lya-Szeg\H{o} inequality \cite{PS}. The following one is Corollary $3$ on Page $66$ in \cite{DM}.  \\ 

\noindent {\bf Lemma $3$} \cite{DM}. Let the real numbers $a_k$ and $b_k$ ($k = 1, 2, \cdots, s$) satisfy 
$0 < m_1 \leq a_k \leq M_1$ and $0 < m_2 \leq b_k \leq M_2$. Then
$$\sum_{k = 1}^s a_k^2 \, \sum_{k = 1}^s b_k^2 \leq \frac{(M_1M_2 + m_1m_2)^2}{4m_1m_2M_1M_2}\left(\sum_{k = 1}^s a_kb_k \right)^2.$$
If $M_1M_2 > m_1m_2$, then the equality sign holds in above inequality if and only if 
$$\nu = \frac{M_1m_2}{M_1m_2 + m_1M_2} s$$
is an integer; while, at the same time, for $\nu$ values of $k$ one has $(a_k, b_k) = (m_1, M_2)$ and for the remaining $s - \nu$ values of 
$k$ one has $(a_k, b_k) = (M_1, m_2)$. If $M_1M_2 = m_1m_2$, the equality always holds.  \\

Lemma $4$ below is from \cite{M}.\\

\noindent {\bf Lemma $4$} \cite{M}. Let $G$ be a balanced bipartite graph of order $2n$ with bipartition ($A$, $B$). If
$d(x) + d(y) \geq n + 1$ for any $x \in A$ and any $y \in B$ with $xy \not \in E$, then $G$ is Hamiltonian. \\

Lemma $5$ below is from \cite{J}. \\ 

\noindent {\bf Lemma $5$} \cite{J}. Let $G$ be a $2$-connected bipartite graph with bipartition ($A$, $B$), where $|A| \geq |B|$. If
each vertex in $A$ has degree at least $s$ and each vertex in $B$ has degree at least $t$, then $G$ contains a cycle
of length at least $2 \min (|B|, s + t - 1, 2s - 2)$. \\

\noindent {\bf 3 \, Proofs } \\
 
\noindent {\bf Proof of Theorem $1$.} Let $G$ be a $k$-connected ($k \geq 2$) graph with $n \geq 3$ vertices and $e$ edges satisfying the conditions 
in Theorem $1$. Suppose $G$ is not Hamiltonian.
Then Lemma $1$ implies that $\beta \geq k + 1$. Also, we have that $n \geq 2 \delta + 1 \geq 2 k + 1$ otherwise 
$\delta \geq k \geq n/2$
and $G$ is Hamiltonian.
 Let $I_1 := \{\, u_1, u_2, ..., u_{\beta} \,\}$ be a maximum independent set in $G$. Then
  $I := \{\, u_1, u_2, ..., u_{k + 1} \,\}$ is an independent set in $G$.
  Thus
$$ \sum_{u \in I} d(u) = |E(I, V - I)| \leq \sum_{v \in V - I} d(v).$$ 
Since $\sum_{u \in I} d(u) +  \sum_{v \in V - I} d(v) = 2e$, we have that 
$$\sum_{u \in I} d(u) \leq e \leq \sum_{v \in V - I} d(v).$$
Notice that $0 < \delta \leq d(u) \leq n - k - 1$ for each $u \in I$. Applying Lemma $3$ with $s =  k + 1$, 
$a_i = 1$ and $b_i = d(u_i)$ with $i = 1, 2, ... , (k + 1)$,  $m_1 = 1 > 0$, $M_1 = 1$,
$m_2 = \delta > 0$, and $M_2 = n - k - 1$, we have 
$$\sum_{i = 1}^{k + 1} 1^2 \, \sum_{i = 1}^{k + 1} d^2(u_i) \leq \frac{(\delta + n - k - 1)^2}{4 \delta (n - k - 1)}\left(\sum_{i = 1}^{k + 1}d(u_i)\right)^2 \leq \frac{(e (\delta + n - k - 1))^2}{4 \delta (n - k - 1)}.$$ 
\noindent Thus $$\sum_{i = 1}^{k + 1} d^2(u_i) \leq \frac{(e (\delta + n - k - 1))^2}{4 \delta (n - k - 1)(k + 1)}.$$ 
\noindent Therefore $$(n - k - 1)\Delta^2 + \frac{(e (\delta + n - k - 1))^2}{4 \delta (n - k - 1)(k + 1)}$$ 
$$\leq M_1 =  \sum_{v \in V - I} d^2(v) + \sum_{u \in I} d^2(u)$$
$$\leq (n - k - 1)\Delta^2 + \frac{(e (\delta + n - k - 1))^2}{4 \delta (n - k - 1)(k + 1)}.$$
\noindent Hence $d(v) = \Delta$ for each $v \in V - I$,
$$\sum_{i = 1}^{k + 1} 1^2 \, \sum_{i = 1}^{k + 1} d^2(u_i) = \frac{(\delta + n - k - 1)^2}{4 \delta (n - k - 1)}\left(\sum_{i = 1}^{k + 1}d(u_i)\right)^2,$$
and $\sum_{i = 1}^{k + 1}d(u_i) = e$ which implies $\sum_{v \in V - I}d(v) = e$ and $G$ is a bipartite graph with partition sets of $I$ and $V - I$. \\

The remaining proofs are divided into two cases. \\

{\bf Case 1.} $\delta =  n - k - 1$. \\

In this case, we have $d(u) = \delta$ for each $u$ in $I$ and thereby $\delta (k + 1) = |E(I, V - I)| = \Delta (n - k - 1) \geq \delta (n - k - 1) $. 
Thus
$n \leq 2 k + 2$. Since $n \geq 2 k + 1$, we have $n = 2 k + 2$ or $n = 2 k + 1$. If $n = 2 k + 2$, then Lemma $4$ implies that $G$ is Hamiltonian,
a contradiction. If $n = 2 k + 1$, then $G$ is $K_{k, \, k + 1}$. \\

{\bf Case 2.} $\delta <  n - k - 1$. \\

In this case, Set $P = \{\, x : x \in I, d(x) = n - k - 1 \,\}$ and $Q = \{\, y : y \in I, d(y) = \delta \,\}$. From Lemma $3$, we have
$|P| = \frac{\delta (k + 1)}{\delta + n - k - 1}$, 
$|Q| = (k + 1) - |P| = \frac{(n - k - 1)(k + 1)}{\delta + n - k - 1}$, and
$I = P \cup Q$.  Choose one vertex $x$ in $P$ 
and one vertex $z$ in $V - I$. Then
$n - k - 1 = d(x) \leq \Delta = d(z) \leq k + 1$. Thus $n \leq 2 k + 2$. Since $n \geq 2 k + 1$, we have $n = 2 k + 2$ or $n = 2 k + 1$. 
If $n = 2 k + 2$, then Lemma $4$ implies that $G$ is Hamiltonian,
a contradiction. If $n = 2 k + 1$, then $G$ is $K_{k, \, k + 1}$ which implies $n - k - 1 = \delta$, a contradiction. \\ 

This completes the proof of Theorem $1$. \\

The proof of Theorem $2$ is similar to the proof of Theorem $1$. For the sake of completeness, 
we still present a full proof of Theorem $2$ below. \\

\noindent {\bf Proof of Theorem $2$.} Let $G$ be a $k$-connected ($k \geq 1$) graph with $n \geq 9$ vertices and $e$ edges satisfying the conditions 
in Theorem $2$. Suppose $G$ is not traceable.
Then Lemma $2$ implies that $\beta \geq k + 2$. Also, we have that $n \geq 2 \delta + 2 \geq 2 k + 2$ otherwise 
$\delta \geq k \geq (n - 1)/2$
and $G$ is traceable.
 Let $I_1 := \{\, u_1, u_2, ..., u_{\beta} \,\}$ be a maximum independent set in $G$. Then
  $I := \{\, u_1, u_2, ..., u_{k + 2} \,\}$ is an independent set in $G$.
  Thus
$$ \sum_{u \in I} d(u) = |E(I, V - I)| \leq \sum_{v \in V - I} d(v).$$ 
Since $\sum_{u \in I} d(u) +  \sum_{v \in V - I} d(v) = 2e$, we have that 
$$\sum_{u \in I} d(u) \leq e \leq \sum_{v \in V - I} d(v).$$
Notice that $0 < \delta \leq d(u) \leq n - k - 2$ for each $u \in I$. Applying Lemma $3$ with $s =  k + 2$, 
$a_i = 1$ and $b_i = d(u_i)$ with $i = 1, 2, ... , (k + 2)$,  $m_1 = 1 > 0$, $M_1 = 1$,
$m_2 = \delta > 0$, and $M_2 = n - k - 2$, we have 
$$\sum_{i = 1}^{k + 2} 1^2 \, \sum_{i = 1}^{k + 2} d^2(u_i) \leq \frac{(\delta + n - k - 2)^2}{4 \delta (n - k - 2)}\left(\sum_{i = 1}^{k + 2}d(u_i)\right)^2 \leq \frac{(e (\delta + n - k - 2))^2}{4 \delta (n - k - 2)}.$$ 
\noindent Thus $$\sum_{i = 1}^{k + 2} d^2(u_i) \leq \frac{(e (\delta + n - k - 2))^2}{4 \delta (n - k - 2)(k + 2)}.$$ 
\noindent Therefore $$(n - k - 2)\Delta^2 + \frac{(e (\delta + n - k - 2))^2}{4 \delta (n - k - 2)(k + 2)}$$ 
$$\leq M_1 =  \sum_{v \in V - I} d^2(v) + \sum_{u \in I} d^2(u)$$
$$\leq (n - k - 2)\Delta^2 + \frac{(e (\delta + n - k - 2))^2}{4 \delta (n - k - 2)(k + 2)}.$$
\noindent Hence $d(v) = \Delta$ for each $v \in V - I$,
$$\sum_{i = 1}^{k + 2} 1^2 \, \sum_{i = 1}^{k + 2} d^2(u_i) = \frac{(\delta + n - k - 2)^2}{4 \delta (n - k - 2)}\left(\sum_{i = 1}^{k + 2}d(u_i)\right)^2,$$
and $\sum_{i = 1}^{k + 2}d(u_i) = e$ which implies $\sum_{v \in V - I}d(v) = e$ and $G$ 
is a bipartite graph with partition sets of $I$ and $V - I$. \\ 

The remaining proofs are divided into two cases. \\

{\bf Case 1.} $\delta =  n - k - 2$. \\

In this case, we have $d(u) = \delta$ for each $u$ in $I$ and thereby $\delta (k + 2) = |E(I, V - I)| = \Delta (n - k - 2)  \geq \delta (n - k - 2)$. 
Thus
$n \leq 2 k + 4$. Since $n \geq 2 k + 2$, we have $n = 2 k + 4$, $n = 2 k + 3$, or $n = 2 k + 4$. 
If $n = 2 k + 4$, then $k \geq 3$ since $n \geq 9$. Thus Lemma $4$ implies that $G$ is Hamiltonian and thereby $G$ is traceable,
a contradiction. If $n = 2 k + 3$, then $k \geq 3$ since $n \geq 9$. Thus Lemma $5$ implies that $G$ has a cycle of length at least $(n - 1)$ 
and thereby $G$ is traceable,
a contradiction. If $n = 2 k + 2$, then $G$ is $K_{k, \, k + 2}$. \\

{\bf Case 2.} $\delta <  n - k - 2$. \\

In this case, Set $P = \{\, x : x \in I, d(x) = n - k - 2 \,\}$ and $Q = \{\, y : y \in I, d(y) = \delta \,\}$. From Lemma $3$, we have 
$|P| = \frac{\delta (k + 2)}{\delta + n - k - 2}$, $|Q| = (k + 2) - |P| = \frac{(n - k - 2)(k + 2)}{\delta + n - k - 2}$, and 
$I = P \cup Q$. Choose one vertex $x$ in $P$ 
and one vertex $z$ in $V - I$. Then
$n - k - 2 = d(x) \leq \Delta = d(z) \leq k + 2$. Thus $n \leq 2 k + 4$. Since $n \geq 2 k + 2$, we have $n = 2 k + 4$, $n = 2 k + 3$, or $n = 2 k + 2$.
If $n = 2 k + 4$, then $k \geq 3$ since $n \geq 9$. Thus Lemma $4$ implies that $G$ is Hamiltonian and thereby $G$ is traceable,
a contradiction. If $n = 2 k + 3$, then $k \geq 3$ since $n \geq 9$. Thus Lemma $5$ implies that $G$ has a cycle of length at least $(n - 1)$ 
and thereby $G$ is traceable,
a contradiction. If $n = 2 k + 2$, then $G$ is $K_{k, \, k + 2}$ which implies $n - k - 2 = \delta$, a contradiction.  \\

This completes the proof of Theorem $2$. \\

 \noindent {\bf Proof of Theorem $3$.} 
 Let $G$ be a graph with $n$ vertices, $e$ edges, and $\delta \geq 1$. Clearly, $\beta < n$. 
 Let $I := \{\, u_1, u_2, ..., u_{\beta} \,\}$ be a maximum independent set in $G$. Then 
$$ \sum_{u \in I} d(u) = |E(I, V - I)| \leq \sum_{v \in V - I} d(v).$$ 
Since $\sum_{u \in I} d(u) +  \sum_{v \in V - I} d(v) = 2e$, we have that 
$$\sum_{u \in I} d(u) \leq e \leq \sum_{v \in V - I} d(v).$$
Notice that $0 < \delta \leq d(u) \leq n - \beta$ for each $u \in I$. Applying Lemma $3$ with $s =  \beta$, 
$a_i = 1$ and $b_i = d(u_i)$ with $i = 1, 2, ... , \beta$,  $m_1 = 1 > 0$, $M_1 = 1$,
$m_2 = \delta > 0$, and $M_2 = n - \beta$, we have 
$$\sum_{i = 1}^{\beta} 1^2 \, \sum_{i = 1}^{\beta} d^2(u_i) \leq \frac{(\delta + n - \beta)^2}{4 \delta (n - \beta)}\left(\sum_{i = 1}^{\beta}d(u_i)\right)^2 \leq \frac{(e (\delta + n - \beta))^2}{4 \delta (n - \beta)}.$$ 
\noindent Thus $$\sum_{i = 1}^{\beta} d^2(u_i) \leq \frac{(e (\delta + n - \beta))^2}{4 \delta (n - \beta)\beta}.$$ 
\noindent Therefore  
$$M_1 =  \sum_{v \in V - I} d^2(v) + \sum_{u \in I} d^2(u)$$
$$\leq (n - \beta)\Delta^2 + \frac{(e (\delta + n - \beta))^2}{4 \delta (n - \beta)\beta}.$$
\noindent If $$M_1 = (n - \beta)\Delta^2 + \frac{(e (\delta + n - \beta))^2}{4 \delta (n - \beta)\beta},$$
\noindent then $d(v) = \Delta$ for each $v \in V - I$,
$$\sum_{i = 1}^{\beta} 1^2 \, \sum_{i = 1}^{\beta} d^2(u_i) = \frac{(\delta + n - \beta)^2}{4 \delta (n - \beta)}\left(\sum_{i = 1}^{\beta}d(u_i)\right)^2,$$
and $\sum_{i = 1}^{\beta}d(u_i) = e$ which implies $\sum_{v \in V - I}d(v) = e$ and $G$ is a bipartite graph with partition sets of $I$ and $V - I$. \\
 
The remaining proofs are divided into two cases. \\

{\bf Case 1.} $\delta =  n - \beta$. \\

In this case, we have $d(u) = \delta$ for each $u$ in $I$ and thereby $G$ is $K_{\beta, \, n - \beta}$. \\

{\bf Case 2.} $\delta <  n - \beta$. \\

In this case, Set $P = \{\, x : x \in I, d(x) = n - \beta \,\}$ and $Q = \{\, y : y \in I, d(y) = \delta \,\}$. From Lemma $3$, we have 
$I = P \cup Q$, $|P| = \frac{\delta \beta}{\delta + n - \beta}$ which is an integer,
 and $|Q| = \beta - |P| = \frac{(n - \beta)\beta}{\delta + n - \beta}$ which is an integer.   \\

\noindent Suppose $G$ is $K_{\beta, \, n - \beta}$. Since $V - I$ is independent in $G$,  $n - \beta \leq \beta$ and thereby
$\delta = n - \beta$. A simple computation can verify that 
$$M_1 = (n - \beta)\Delta^2 + \beta (n - \beta)^2$$
$$= (n - \beta)\Delta^2 + \frac{(e (\delta + n - \beta))^2}{4 \delta (n - \beta)\beta}.$$

\noindent Suppose $G$ is a bipartite graph with partition sets of $I$ and $V - I$ such that $|I| = \beta$, $\delta <  n - \beta$, $d(v) = \Delta$ for each vertex $v$ in $v - I$, and 
$I = P \cup Q$, where $P = \{\, x : x \in I, d(x) = n - \beta \,\}$,  $Q = \{\, y : y \in I, d(y) = \delta \,\}$,  
$|P| = \frac{\delta \beta}{\delta + n - \beta}$ is an integer, and $|Q| = \frac{(n - \beta)\beta}{\delta + n - \beta}$ is an integer. Then
$$e = |P| (n - \beta) + |Q| \delta = \frac{\delta \beta (n - \beta)}{\delta + n - \beta} + \frac{\delta \beta (n - \beta)}{\delta + n - \beta}  
= \frac{2\delta \beta (n - \beta)}{\delta + n - \beta}.$$
\noindent Thus   
$$M_1 = \sum_{v \in V - I} d^2(v) + \sum_{u \in I} d^2(u) = (n - \beta)\Delta^2 + |P| (n - \beta)^2 + |Q| \delta^2 $$ 
$$= (n - \beta)\Delta^2 + \delta (n - \beta) \beta = (n - \beta)\Delta^2 + \frac{(e(\delta + n - \beta))^2}{4 \delta (n - \beta) \beta}.$$

This completes the proof of Theorem $3$. \\

\end{document}